**Invited article**

# SOM-based algorithms for qualitative variables

**Marie Cottrell, Smaïl Ibbou, Patrick Letrémy**

SAMOS-MATISSE UMR 8595

90, rue de Tolbiac, F-75634 Paris Cedex 13, France

Telephone and fax number: +33 1 44 07 89 22

cottrell@univ-paris1.fr

*Abstract*: It is well known that the SOM algorithm achieves a clustering of data which can be interpreted as an extension of Principal Component Analysis, because of its topology-preserving property. But the SOM algorithm can only process real-valued data. In previous papers, we have proposed several methods based on the SOM algorithm to analyze categorical data, which is the case in survey data. In this paper, we present these methods in a unified manner. The first one (Kohonen Multiple Correspondence Analysis, KMCA) deals only with the modalities, while the two others (Kohonen Multiple Correspondence Analysis with individuals, KMCA_ind, Kohonen algorithm on DISJonctive table, KDISJ) can take into account the individuals, and the modalities simultaneously.

## 1. Introduction to data analysis

### *1.1 Data Analysis, classical methods*

The extraction of knowledge from large databases is an essential part of data analysis and data mining. The original information has to be summed up, by simplifying the amount of rough data, computing some basic features, and giving visual representation.

To analyze, sum up, and represent multidimensional data involving quantitative (continuous) variables and qualitative (nominal, ordinal) variables, experts have numerous performing methods at disposal, experimented and developed in most statistical software packages. Data analysis consists in building simplified representations, so as to highlight the relations, the main characteristics, and the internal structure of the data set.

One can distinguish two main groups of classical techniques: factorial methods and classification methods.

The factorial methods are inherently linear: they consist in looking for vectorial sub-spaces and new axes, in reducing the data dimension, while keeping the greater part of the information. The Principal Component Analysis (PCA) is a very popular technique which projects quantitative data on significant axes, and the Correspondence Analysis provides an analysis of the relations which can exist between the modalities of all the qualitative variables by completing a simultaneous projection of the modalities. There are two variants, the Factorial Correspondence Analysis (FCA) for two variables and the Multiple Correspondence Analysis (MCA) for more than two variables.

On the other hand, the classification methods are numerous and diverse. They group the data into homogeneous clusters. The most employed ones are the Ascending Hierarchical Classification, in which the number of classes is not fixed *a priori*, and the Moving Centers Method (also called k-Means Method) which gathers the data into an *a priori* given number of classes.

These classical methods are presented for example in Hotelling (1933), Cooley and Lohnes (1971), Morrison (1976), Mardia et al. (1979), as to the Principal Component Analysis, and in Burt (1950), Benzecri (1973), Greenacre (1984), Lebart et al. (1984) as to the Correspondence Analysis. The classification methods can be found in Anderberg (1973), Hartigan (1975). In French, Lebart et al. (1995) or Saporta (1990) present these methods with a lot of examples.

### *1.2 Neural Methods*

Recently, - from the 80s - , new methods have appeared, known as neural methods. They come from interdisciplinary studies which involve biologists, physicists, computer scientists, signal processing experts, psychologists, mathematicians, and statisticians. See for example Rumelhart and McClelland (1986), Hertz et al. (1991), Haykin (1994) for foundations of Neural Networks.

These methods have rapidly encountered a large success, in particular because they originally appear as “black boxes” able to do anything in numerous application fields. After having gone beyond the excessive enthusiasm these methods first aroused and overcome some difficulties of implementation, researchers and users now have many toolboxes of alternative techniques, generally non linear and iterative.

For example, the reader can consult the book by Ripley (1996) which integrates the neural techniques into the statistical methods. Many statisticians include these methods into their tools. See for example the neural packages of the main statistical software, like SAS, SYLAB, STATLAB, S+ and of the computing software (MATLAB, R, GAUSS, etc.).

The main point is that when the standard linear statistical methods are not appropriate, due to the intrinsic structure of the observations, the neural models which are highly non linear can be very helpful.

A very popular neural model is the Multilayer Perceptron (Rumelhart and McClelland (1986), Hertz et al. (1991). One of its nice properties is that it admits both quantitative and qualitative variables as inputs. It works very well to achieve classification, discrimination, non linear regression analysis, short term forecasting and function approximation. But it uses supervised learning, which means that the labels of the classes for each datum, i.e. the desired outputs, have to be known *a priori*. On the other hand, it is not very easy to interpret the model, to find the most relevant variables, and to propose a suitable typology of the observations from the results. Even if many advances have been made since the « black box » epoch, to better choose the architecture of the network, to prune the non-significant connections, to offer the extraction of rules, and so on, this kind of models is not very appropriate for the interpretation, representation and visualization of the data, for which there is no *a priori* knowledge.

So in data analysis, unsupervised methods are very attractive and in particular the Kohonen algorithm is nowadays widely used in this framework (Kohonen, 1984, 1993, 1995, 1997). It achieves both tasks of “projection” and classification. Let us recall the definition of this algorithm, also called SOM (Self-Organizing Map). In its genuine form, it processes quantitative real-valued data, in which each observation is described by a real vector. For example the quantitative variables can be ratios, quantities, measures, indices, coded by real numbers. For the moment, we do not consider the qualitative variables which can be present in the database.

In this case (only quantitative variables), we consider a set of $N$ observations, in which each individual is described by $p$ quantitative real-valued variables. The main tool is a Kohonen network, generally a two-dimensional grid, with $n$ by $n$ units, or a one-dimensional string with $n$ units. The data are arranged in a table $X$ with $N$ rows and $p$ columns. The rows of table $X$ are the inputs of the SOM algorithm. After learning, each unit $u$ is represented in the $R^p$ space by its weight vector $C_u$ (or *code vector*). Then each observation is classified by a nearest neighbor method: observation $i$ belongs to class $u$ if and only if the code vector $C_u$ is the closest among all the code vectors. The distance in $R^p$ is the Euclidean distance in general, but one can choose a different one depending on the application.

Compared to any other classification method, the main characteristic of the Kohonen classification is the conservation of the topology: after learning, « close » observations are associated to the same class or to « close » classes according to the definition of the neighborhood in the Kohonen network.

There are a lot of applications of this algorithm to real data, see for example the immense bibliography available on the WEB page at http:/ www.cis.hut.fi/.

All these studies show that the SOM algorithm has numerous nice properties : the representation on the grid or on the string is easy to interpret, the topology conservation gives some “order” among the classes, it is possible to use data with missing values and the classification algorithm is quick and efficient (De Bodt et al. 2003).

### *1.3 Kohonen map versus PCA*

The Kohonen map built with the rows of a data table *X* can be compared to the linear projections achieved by a Principal Component Analysis (PCA) on the successive principal axes. See for example Blayo and Desmartines (1991, 1992) or Kaski (1997). To find these axes, one has to compute the eigenvalues and the eigenvectors of the matrix *X’X* where *X’* is the transpose matrix of *X*. Each eigenvalue is equal to the part of the total inertia represented on the corresponding axis. The axes are ordered according to the decreasing eigenvalues. So the first two axes correspond to the two largest eigenvalues, and provide the best projection. However it is often necessary to take into account several 2-dimensional PCA projections to get a good representation of the data, while there is only one Kohonen map.

**Let us emphasize this point: if *X* is the data matrix, the PCA is achieved by diagonalizing the matrix *X'X*, while the Kohonen map is built with the rows of the data matrix *X*.**

## 2. Preliminary comments about qualitative variables

### *2.1 Qualitative variables*

In real-world applications, the individuals can also be described by variables which have a qualitative (or categorical) nature. This is for example the case when the data are collected after a survey, in which people have to answer a number of questions, each of them having a finite number of possible modalities (i.e. sex, professional group, level of income, kind of employment, place of housing, type of car, level of education, etc.). The values of these variables are nominal values.

Sometimes the modalities of these variables are encoded by numerical values 1, 2, 3,... that could be viewed as numerical values, but it is well known that this is not adequate in general. The encoding values are not always comparable and/or the codes are neither necessarily ordered nor regularly spaced (for example, is blue color less than brown color? how to order the types of car, the places of housing, ...?). Most of the time, using the encoding of the modalities as quantitative (numerical) variables has no meaning. Even if the encoding values correspond to an increasing or decreasing progression, (from the poorest to the richest, from the smallest to the largest, etc.), it is not correct to use these values as real values, except if a linear scaling is used, that is if modality 2 is exactly halfway between modalities 1 and 3, etc. *So the qualitative data need a specific treatment*. When the data include both types of variables, the first idea is to achieve the classification using only the quantitative variables and then, to cross the classification with the other qualitative variables.

### *2.2 Crossing a classification with a qualitative variable*

Let us suppose that the individuals have been clustered into classes on the basis of the real-valued variables which describe them, using a Kohonen algorithm.

We can consider some non-used qualitative variables to answer some questions: what is the nature of the observations in a given class, is there a characteristic common to the neighboring classes, can we qualify a group of classes by a qualitative variable? A first method can be to extract the observations of a given class and analyze them with statistical software, by computing means and variances, for the quantitative variables (used for the classification) and frequencies for the qualitative variables. That gives an answer to the first question, but we lose the neighborhood properties of the Kohonen map.

In order to complete the description, we can study the repartition of each qualitative variable within each class. Let be $Q$ a qualitative variable with $m$ modalities. In each cell of the Kohonen map, we draw a frequency pie, where each modality is represented by a grey level occupying an area proportional to its frequency in the corresponding class. See in Fig. 1 an example of a frequency pie.

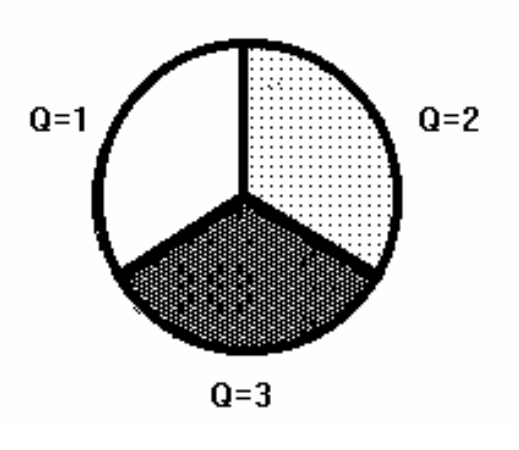

Fig. 1: *A frequency pie, when there are 3 modalities with frequencies 1/3.*

So, in this way, by representing the frequencies of each modality across the map, we can add some information which contributes to the description of the class, while keeping the topology-conservation property of the map and making clear the continuity as well as the breakings.

### *2.3. Example I: the country database, POP_96.*

Let us present a classical example. The POP_96 database contains seven ratios measured in 1996 on the macroeconomic situation of 96 countries: annual population growth, mortality rate, analphabetism rate, population proportion in high school, GDP per head, unemployment rate and inflation rate. This dataset was first used in Blayo and Demartines (1991, 1992) in the context of data analysis by SOMs. The 96 countries which are described by these 7 real-valued variables, are first classified on a Kohonen map with 6 by 6 units.

The initial choice of the number of classes is arbitrary, and there does not exist any method to choose the size of the grid in a perfect way, even if some authors suggest using growing or shrinking approaches (Blackmore and Mikkulainen, 1993, Fritzke, 1994, Koikkalainen and Oja, 1990). To exploit the stochastic features of the SOM algorithm, and to obtain a good clustering and a good organization, it is proven to be more efficient to deal with “large” maps. But we can guess that the « relevant » number of classes would often be smaller than the size $n \times n$ of the grid. On the other hand, it is neither easy nor useful to give interpretation and description of a too large number of classes. So we propose (Cottrell et al., 1997) to reduce the number of classes by means of an Ascending Hierarchical Classification of the code vectors using the Ward distance for example (Lance and Williams, 1967, Anderberg, 1973).

In this way, we define two embedded classifications, and can distinguish the classes (Kohonen classes or « micro-classes ») and the « macro-classes » which group together some of the « micro-classes ». To make this two-level classification visible, we assign to each « macro-class » some color or grey level.

The advantage of this double classification is the possibility to analyze the data set at a « macro » level where general features emerge and at a « micro » level to determine the characteristics of more precise phenomena and especially the paths to go from one class to another one.

In the applications that we treated, the « macros-classes » always create connected areas in the grid. This remark is very attractive because it confirms the topological properties of the Kohonen maps.

See in Fig.2, a representation of the « micro-classes » grouped together to constitute 7 « macro-classes ». The macro-class in the top left-hand corner groups the OECD, rich and very developed countries, the very poor ones are displayed on the right, ex-socialist countries are not very far from the richest, the very inflationist countries are in the middle at the top, and so on.

| | | | | | |
|---|---|---|---|---|---|
| Japan<br>Switzerland | USA<br>Luxembourg | Cyprus<br>South Korea | Russia<br>Ukraine | Brazil | Angola |
| Germany<br>Belgium<br>Denmark<br>France<br>Norway<br>Netherlands<br>Sweden | Australia<br>Canada<br>Iceland | U Arab<br>Emirates<br>Israel<br>Singapore | | | Afghanistan<br>Mozambique<br>Yemen |
| Spain<br>Finland<br>Ireland<br>Italy<br>New-Zealand<br>United Kingdom | Greece<br>Malta<br>Portugal | Argentina<br>Bahrain<br>Philippines | Mongolia<br>Peru | Swaziland | Mauritania<br>Sudan |
| Croatia<br>Hungary<br>Romania<br>Slovenia<br>R. Czech | Bulgaria<br>Poland<br>Yugoslavia | Jamaica<br>Sri Lanka | Tunisia<br>Turkey | Saudi Arabia<br>Bolivia<br>El Salvador<br>Syria | Comoros<br>Ivory Coast<br>Ghana<br>Morocco<br>Pakistan |
| Moldavia<br>Uruguay | Chili<br>Colombia | Albania<br>Panama | | Algeria<br>Egypt<br>Nicaragua | Cameroon<br>Laos<br>Nigeria |
| China<br>Fiji<br>Malaysia<br>Mexico<br>Thailand | Costa Rica<br>Ecuador<br>Guyana<br>Paraguay<br>Venezuela | Vietnam | South Africa<br>Lebanon<br>Macedonia | Indonesia<br>Iran<br>Namibia<br>Zimbabwe | Haiti<br>Kenya |

Fig. 2: *The 36 Kohonen classes, grouped into 7 macro-classes,* 600 *iterations*.

In order to describe the countries more precisely, we can for example consider an extra qualitative variable, classically defined by the economists. The IHD variable is the Index of Human Development which takes into account a lot of data which can qualify the way of life, the cultural aspects, the security, the number of physicians, of theaters, etc. So the qualitative variable is the IHD index, with 3 levels: low (yellow), medium (green) and high (blue). See Fig. 3 the repartition of this extra qualitative variable across the Kohonen map.

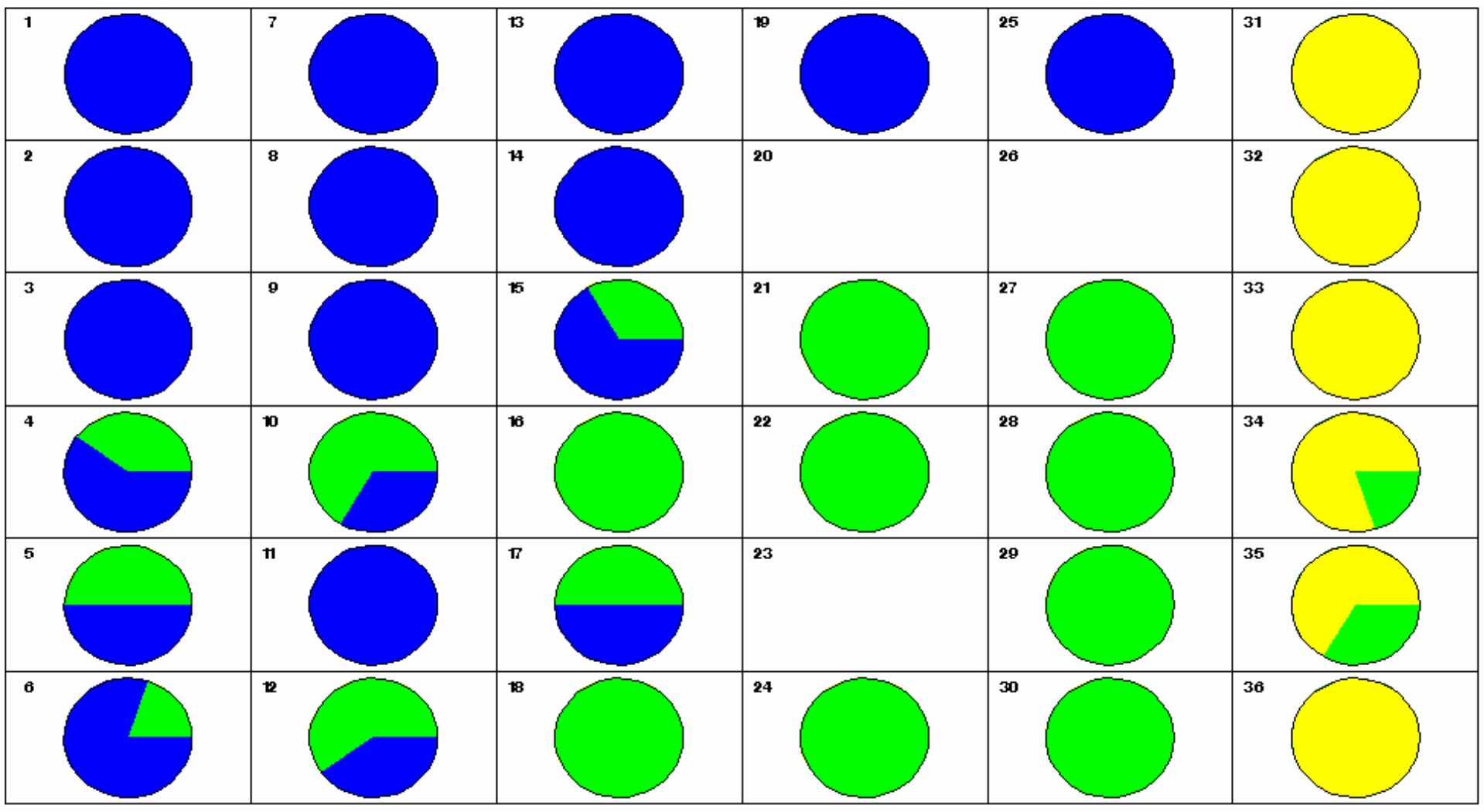

Fig. 3: *In each cell the frequency pie of the variable IHD is represented. We can observe that units 1, 7, 13, 19, 25, 2, 8, 14, 3, 9, 15, 4, 5, 11, 17, 6 are mainly countries with high levels of IHD. Units 31, 32, 33, 34, 35, 36 have the lowest level. This is fully coherent with the above classification (Fig. 2).*

In sections 3 and 4, we introduce the categorical database and the classical method to deal with them. In sections 5, we define the algorithm KMCA, which is a Kohonen-based method. It achieves a first classification which is easy to visualize, a nice representation of the modalities, which can be clustered into macro-classes. This algorithm is applied to a real-world database in section 6. In sections 7 and 8, we define two algorithms which take into account the individuals together with the modalities. The same example is used along all these sections. Sections 9 and 10 present two other examples, one of them is a toy example, while the other is taken from real data. Section 12 is devoted to conclusion and perspectives.

## 3. Database with only qualitative variables

From now, all the observations are described by $K$ qualitative variables, each of them having a certain number of modalities, as in a survey.

Let us define the data and introduce the basic notations. Let us consider an $N$-sample of individuals and $K$ variables or questions. Each question has $m_k$ possible modalities (or answers, or levels). The individuals answer each question $k$ ($1 \le k \le K$) by choosing only one modality among the $m_k$ modalities. For example, if we assume that $K = 3$ and $m_1 = 3$ , $m_2 = 2$ and $m_3 = 3$, then an answer of an individual could be (0,1,0|0,1|1,0,0), where 1 corresponds to the chosen modality for each question.

So, if $M = \sum_{k=1}^{K} m_k$ is the total number of modalities, each individual is represented by a row $M$-vector with values in {0, 1}. There is only one 1 between the 1st component and the $m_1$-th one, only one 1 between the ($m_1$+1)-th-component and the ($m_1$+$m_2$)-th-one and so on.

To simplify, we can enumerate all the modalities from 1 to $M$ and denote by $Z_j$, ($1 \le j \le M$) the column vector constructed by grouping the $N$ answers to the $j$-th modality. The $i$-th element of the vector $Z_j$ is 1 or 0, according to the choice of the individual $i$ (it is 1 if and only if the individual $i$ has chosen the modality $j$).

Then we can define a ($N\times M$) matrix $D$ as a logical canonical matrix whose columns are the $Z_j$ vectors. It is composed of $K$ blocks where each ($N\times m_k$) block contains the $N$ answers to the question $k$. One has:

$$D = \left(Z_1, \cdots, Z_{m_1}, \cdots, Z_j, \cdots, Z_M\right)$$

The ($N \times M$) data matrix $D$ is called the *complete disjunctive table* and is denoted by

$$D = (d_{ij}),\ i = 1, ..., N,\ j = 1, ..., M.$$

This table $D$ contains all the information about the individuals. This table is the rough result of any survey.

| | $m_1$ | | | $m_2$ | | $m_3$ | | |
|---|---|---|---|---|---|---|---|---|
| *Ind* | 1 | 2 | 3 | 1 | 2 | 1 | 2 | 3 |
| 1 | 0 | 1 | 0 | 0 | 1 | 1 | 0 | 0 |
| 2 | 1 | 0 | 0 | 1 | 0 | 0 | 1 | 0 |
| … | | | | | | | | |
| … | | | | | | | | |
| *i* | 0 | 0 | 1 | 0 | 1 | 0 | 0 | 1 |
| | | | | | | | | |
| | | | | | | | | |
| *N* | 0 | 0 | 1 | 1 | 0 | 0 | 0 | 1 |

Table 1*: Example of a Complete Disjunctive Table.*

If we want to remember who answered what, it is essential use this table, see later in section 5. But if we only have to study the *relations between the K variables (or questions)*, we can sum up the data in a cross-tabulation table, *called Burt matrix*, defined by

$$B = D'D$$

where $D'$ is the transpose matrix of $D$. The matrix $B$ is a ($M \times M$) symmetric matrix and is composed of $K \times K$ blocks, in which the ($k$, $l$) block $B_{kl}$ (for $1 \le k, l \le K$) is the contingency table which crosses the question $k$ and the question $l$. The block $B_{kk}$ is a diagonal matrix, whose diagonal entries are the numbers of individuals who have respectively chosen the modalities 1, ... , $m_k$, for question $k$.

The Burt table can be represented as below. It can be interpreted as a generalized contingency table, when there are more than 2 types of variables to simultaneously study.

From now, we denote the entries of the matrix $B$ by $b_{jl}$, whatever the questions which contain the modalities $j$ or $l$. The entry represents the number of individuals which chose both modalities $j$ and $l$. According to the definition of the data, if $j$ and $l$ are two different modalities of a same question, $b_{jl} = 0$, and if $j = l$ , the entry $b_{jj}$ is the number of individuals who chose modality $j$. In that case, we use only one sub-index and write down $b_j$ instead of $b_{jj}$. This number is nothing else than the sum of the elements of the vector $Z_j$. Each row of the matrix $B$ characterizes *a modality of a question* (also called *variable*). Let us represent below the Burt table for the same case as for the disjunctive table ($K = 3$, $m_1 = 3$, $m_2 = 2$ and $m_3 = 3$)

| | $Z_1$ | $Z_2$ | $Z_3$ | $Z_4$ | $Z_5$ | $Z_6$ | $Z_7$ | $Z_8$ |
|---|---|---|---|---|---|---|---|---|
| $Z_1$ | $b_1$ | 0 | 0 | $b_{14}$ | $b_{15}$ | $b_{16}$ | $b_{17}$ | $b_{18}$ |
| $Z_2$ | 0 | $b_2$ | 0 | $b_{24}$ | $b_{25}$ | $b_{26}$ | $b_{27}$ | $b_{28}$ |
| $Z_3$ | 0 | 0 | $b_3$ | $b_{34}$ | $b_{35}$ | $b_{36}$ | $b_{37}$ | $b_{38}$ |
| $Z_4$ | $b_{41}$ | $b_{42}$ | $b_{43}$ | $b_4$ | 0 | $b_{46}$ | $b_{47}$ | $b_{48}$ |
| $Z_5$ | $b_{51}$ | $b_{52}$ | $b_{53}$ | 0 | $b_5$ | $b_{56}$ | $b_{57}$ | $b_{58}$ |
| $Z_6$ | $b_{61}$ | $b_{62}$ | $b_{63}$ | $b_{64}$ | $b_{65}$ | $b_6$ | 0 | 0 |
| $Z_7$ | $b_{71}$ | $b_{72}$ | $b_{73}$ | $b_{74}$ | $b_{75}$ | 0 | $b_7$ | 0 |
| $Z_8$ | $b_{81}$ | $b_{82}$ | $b_{83}$ | $b_{84}$ | $b_{85}$ | 0 | 0 | $b_8$ |

Table 2: *Example of a Burt Table.*

One can observe that for each row $j$ (or column since $B$ is symmetric), $\sum_{l} b_{jl} = b_{j\cdot} = K\ b_j$ , since this number is repeated in each block of the matrix $B$ and that $\sum_{j} b_j = \sum_{k=1}^{K} \sum_{l=1}^{m_k} b_l = NK$ . So the total sum of all the entries of $B$ is $b = \sum_{j,l} b_{jl} = K \sum_{j} b_j = K^2 N$ .

In the next section, we describe the classical way to study the relations between the modalities of the qualitative variables, that is the Multiple Correspondence Analysis, which is a kind of factorial analysis.

## 4 . Factorial Correspondence Analysis and Multiple Correspondence Analysis

### *4.1. Factorial Correspondence Analysis*

The classical Multiple Correspondence Analysis (MCA) (Burt, 1950, Benzécri, 1973, Greenacre, 1984, Lebart et al., 1984) is a generalization of the Factorial Correspondence Analysis (FCA), which deals with a Contingency Table. Let us consider only two qualitative variables with respectively $I$ and $J$ modalities. The Contingency Table of these two variables is a $I \times J$ matrix, where entry $n_{ij}$ is the number of individuals which share modality $i$ for the first variable (row variable) and modality $j$ for the second one (column variable).

This case is fundamental, since both Complete Disjunctive Table $D$ and Burt Table $B$ can be viewed as Contingency Tables. In fact, $D$ is the contingency table which crosses a "meta-variable" INDIVIDUAL with $N$ values and a "meta-variable" MODALITY with $M$ values. In the same way, $B$ is clearly the contingency table which crosses a “meta-variable” MODALITY having $M$ values with itself.

Let us define a Factorial Correspondence Analysis, applied to some Contingency Table, (Lebart et al., 1984).

One defines successively

- the table $F$ of the relative frequencies, with entry $f_{ij} = \frac{n_{ij}}{n}$, where $n = \sum_{ij} n_{ij}$
- the margins with entry $f_{i\cdot} = \sum_j f_{ij}$ or $f_{\cdot j} = \sum_i f_{ij}$ ,
- the table $P_R$ of the $I$ row profiles which sum to 1, with entry $p_{ij}^R = \frac{f_{ij}}{\sum_j f_{ij}} = \frac{f_{ij}}{f_{i\cdot}}$,
- the table $P_C$ of the $J$ column profiles which sum to 1, with entry $p_{ij}^C = \frac{f_{ij}}{\sum_i f_{ij}} = \frac{f_{ij}}{f_{\cdot j}}$.

These profiles form two sets of points respectively in $R^J$ and in $R^I$. The means of these two sets are respectively denoted by $\bar{i} = \left( f_{\cdot 1}, f_{\cdot 2}, \cdots, f_{\cdot J} \right)$ and $\bar{j} = \left( f_{1\cdot}, f_{2\cdot}, \cdots, f_{I\cdot} \right)$.

As the profiles are in fact conditional probability distributions ($p_{ij}^R$ is the conditional probability that the first variable has value $i$, given that the second one is equal to $j$, same for $p_{ij}^C$), it is usual to consider the $\chi^2$-distance between rows, and between columns. This distance is defined by

$$\chi^2(i,i') = \sum_j \frac{1}{f_{\cdot j}} \left( \frac{f_{ij}}{f_{i\cdot}} - \frac{f_{i'j}}{f_{i'\cdot}} \right)^2 = \sum_j \left( \frac{f_{ij}}{\sqrt{f_{\cdot j}}\, f_{i\cdot}} - \frac{f_{i'j}}{\sqrt{f_{\cdot j}}\, f_{i'\cdot}} \right)^2 \text{ (distance between rows)}$$

and

$$\chi^2(j,j') = \sum_i \frac{1}{f_{i\cdot}}\left(\frac{f_{ij}}{f_{\cdot j}} - \frac{f_{ij'}}{f_{\cdot j'}}\right)^2 = \sum_i \left(\frac{f_{ij}}{\sqrt{f_{i\cdot}}\, f_{\cdot j}} - \frac{f_{ij'}}{\sqrt{f_{i\cdot}}\, f_{\cdot j'}}\right)^2 \text{ (distance between columns)}$$

Note that each row $i$ is weighted by $f_{i\cdot}$ and that each column is weighted by $f_{\cdot j}$.

So it is possible to compute the inertia of both sets of points :

*Inertia* (*row profiles*) = $\sum_i f_{i\cdot}\ \chi^2(i,\bar{i})$ and *Inertia* (*column profiles*) = $\sum_i f_{\cdot j}\ \chi^2(j,\bar{j})$.

It is easy to verify that these two expressions are equal. This inertia is denoted by $\Im$, and can be written:

$$\Im = \sum_{ij} \frac{\left(f_{ij} - f_{i\cdot}f_{\cdot j}\right)^2}{f_{i\cdot}f_{\cdot j}} = \sum_{ij} \frac{f_{ij}^2}{f_{i\cdot}f_{\cdot j}} - 1 .$$

We can underline two important facts:

1) In order to use the Euclidean distance between rows, and between columns instead of the $\chi^2$-distance, and to take into account the weighting of each row by $f_{i\cdot}$ and of each column by $f_{\cdot j}$, it is very convenient to replace the initial values $f_{ij}$ by corrected values $f_{ij}^c$, by putting down

$$f_{ij}^c = \frac{f_{ij}}{\sqrt{f_{i\cdot}}\sqrt{f_{\cdot j}}} .$$

Let us denote by $F^c$ the matrix whose entries are the $f_{ij}^c$.

2) The inertia of both sets of row profiles and column profiles is exactly equal to $\frac{1}{N}T$ , where $T$ is the usual chi-square statistics which is used to test the independence of both row variable and column variable. Statistics $T$ is also a measure of the deviation from independence.

The Factorial Correspondence Analysis (FCA) is merely a double PCA achieved on the rows and on the columns of this corrected data matrix $F^c$. For the row profiles, the eigenvalues and eigenvectors are computed by the diagonalization of the matrix $F^{c\prime}F^c$. For the column profiles, the eigenvalues and eigenvectors are computed by the diagonalization of the transpose matrix $F^cF^{c\prime}$. It is well known that both matrices have the same eigenvalues and that their eigenvectors are strongly related. It is easy to prove that the total inertia $\Im$ is equal to the sum of the eigenvalues of $F^{c\prime}F^c$ or $F^cF^{c\prime}$. So the FCA decomposes the deviation from independence into a sum of decreasing terms associated to the principal axes of both PCAs sorted out according to the decreasing order of the eigenvalues.

For this FCA, which deals with only two variables, the coupling between the two PCAs is ensured, because they act on two transpose matrices. It is thus possible to simultaneously represent the modalities of both variables.

**According to section 1.3, the diagonalization of the data matrix $F^{c\prime}F^c$ can be approximately replaced by a SOM algorithm in which the row profiles are used as inputs, whereas the diagonalization of $F^cF^{c\prime}$ can be replaced by a SOM algorithm in which the column profiles are used as inputs.**

**This is the key point for defining the SOM algorithms adapted to qualitative variables.**

*4.2. Multiple Correspondence Analysis*

Let us now recall how the classical Multiple Correspondence Algorithm is defined.

**1) In the case of a MCA, when we are interested in the modalities only,** the data table is the Burt table, considered as a contingency table. As explained just before, we consider the c*orrected Burt table* $B^c$, with

$$b_{jl}^c = \frac{b_{jl}}{\sqrt{b_{j\cdot}}\sqrt{b_{\cdot l}}} = \frac{b_{jl}}{K\sqrt{b_j}\sqrt{b_l}},$$

since $b_{j\cdot} = Kb_j$ and $b_{\cdot l} = Kb_l$. As matrices $B$ and $B^c$ are symmetric, the diagonalizations of $B^{c\prime}B^c$ or $B^cB^{c\prime}$ are identical. Then the principal axes of the usual Principal Component Analysis of $B^c$, are the principal axes of the Multiple Correspondence Analysis. It provides a simultaneous representation of the $M$ row vectors, i.e. of the modalities, on several 2-dimensional spaces which give information about the relations between the $K$ variables.

**2) If we are interested in the individuals**, it is necessary use the Complete Disjunctive table $D$, considered as a contingency table (see above).

Let us denote by $D^c$ the corrected matrix, where the entry $d_{ij}^c$ is given by

$$d_{ij}^c = \frac{d_{ij}}{\sqrt{d_{i.}}\sqrt{d_{.j}}} = \frac{d_{ij}}{\sqrt{K}\sqrt{b_j}}.$$

In this case, this matrix is no longer symmetric. The diagonalization of $D^{c}{}'D^{c}$ will provide a representation of the individuals, while the diagonalization of $D^{c}D^{c}{}'$ will provide a representation of the modalities. Both representations can be superposed, and provide the simultaneous representations of individuals and modalities. In that case, it is possible to compute the coordinates of the modalities and individuals.

Let us conclude this short definition of the classical Multiple Correspondence Analysis, by some remarks. MCA is a linear projection method, and provides several two-dimensional maps, each of them representing a small percentage of the global inertia. It is thus necessary to look at several maps at once, the modalities are more or less well represented, and it is not always easy to deduce pertinent conclusions about the proximity between modalities. Related modalities are projected onto neighboring points, but it is possible that neighboring points do not correspond to neighboring modalities, because of the distortion due to the linear projection. Their main property is that each modality is drawn as an approximate center of gravity of the modalities which are correlated with it and of the individuals who share it (if the individuals are available). But the approximation can be very poor and the graphs are not always easy to interpret, as we will see in the examples.

### *4.3. From Multiple Correspondence Analysis to SOM algorithm*

Therefore it is easy to define the new SOM-based algorithms, following section 1.3.

1) In case we want to deal only with the modalities, as the Burt matrix is symmetric, it is sufficient to use a SOM algorithm over the rows (or over the columns) of $B^c$ to achieve a nice representation of all the modalities on a Kohonen map. This remark founds the definition of the algorithm KMCA (Kohonen Multiple Correspondence Analysis), see section 5.

2) If we want to keep the individuals, we can apply the SOM algorithm to the rows of $D^c$, but we will get a Kohonen map for the individuals only. To simultaneously represent the modalities, it is necessary to use some other trick.

Two techniques are defined:

a) KMCA_ind (Kohonen Multiple Correspondence Analysis with Individuals): the modalities are assigned to the classes after training, as supplementary data (see section 7).

b) KDISJ (Kohonen algorithm on DISJunctive table): two SOM algorithms are used on the rows (individuals) and on the columns (modalities) of $D^c$, while being compelled to be associated all along the training (see section 8).

## 5 Kohonen-based Analysis of a Burt Table : algorithm KMCA

In this section, we only take into account the modalities and define a Kohonen-based algorithm, which is analogous to the classical MCA on the Burt Table.

This algorithm was introduced in Cottrell, Letrémy, Roy (1993), Ibbou, Cottrell (1995), Cottrell, Rousset (1997), the PHD thesis of Smaïl Ibbou (1998), Cottrell et al. (1999), Letrémy, Cottrell, (2003). See the references for first presentations and applications.

The data matrix is the corrected Burt Table $B^c$ as defined in Section 4.2. Consider an $n \times n$ Kohonen network (bi-dimensional grid), with a usual topology.

Each unit $u$ is represented by a code vector $C_u$ in $R^{\boldsymbol{M}}$; the code vectors are initialized at random. The training at each step consists of

- presenting at random an input $r(j)$ i.e. a row of the *corrected matrix* $B^c$,
- looking for the winning unit $u_0$, i.e. that which minimizes $\|r(j) - C_u\|^2$ for all units $u$,
- updating the weights of the winning unit and its neighbors by

$$C_u^{new} - C_u^{old} = \varepsilon\, \sigma(u,u_0)\, (r(j) - C_u^{old}).$$

where $\varepsilon$ is the adaptation parameter (positive, decreasing with time), and $\sigma$ is the neighborhood function, with $\sigma(u, u_0)$= 1 if $u$ and $u_0$ are neighbor in the Kohonen network, and = 0 if not. The radius of the neighborhood is also decreasing with time. [1]

After training, each row profile $r(j)$ is represented by its corresponding winning unit. Because of the topology-preserving property of the Kohonen algorithm, the representation of the $M$ inputs on the $n \times n$ grid highlights the *proximity* between the modalities of the $K$ variables.

[1] The adaptation parameter is defined as a decreasing function of the time $t$, which depends on the number $n \times n$ of units in the network, $\varepsilon = \varepsilon_0 / (1 + c_0\, t / n \times n)$. The radius of neighborhood is also a decreasing function of $t$, depending on $n$ and on $T_{max}$ (the total number of iterations) , $\rho(t) = \text{Integer}\left(\frac{n/2}{1+(2n-4)/T_{\max}}\right)$

After convergence, we get an organized classification of all the modalities, where related modalities belong to the same class or to neighboring classes.

We call this method the **Kohonen Multiple Correspondence Analysis (KMCA)** which provides a very interesting alternative to classical Multiple Correspondence Analysis.

## 6. Example I: the country database with qualitative variables

Let us consider the POP_96 database introduced in section 2.3. Now, we consider the 7 variables as qualitative ones, by discretization into classes, and we add the eighth variable IHD as before. The 8 variables are defined as represented in table 3.

| **Heading** | **Modalities** | **Name** |
|---|---|---|
| **Annual population growth** | [-1, 1[, [1, 2[, [2, 3[, ≥3 | ANPR1, ANPR2, ANPR3, ANPR4 |
| **Mortality rate** | [4, 10[10, 40[40, 70[, [70, 100[,≥ 100 | MORT1, MORT2, MORT3, MORT4, MORT5 |
| **Analphabetism rate** | [0, 6[, [6, 20[, [20, 35[, [35, 50[, ≥ 50 | ANALR1, ANALR2, ANALR3, ANALR4, ANALR5 |
| **High school** | ≥ 80, [40, 80[, [4, 40[ | SCHO1, SCHO2, SCHO3 |
| **GDPH** | ≥ 10000, [3000,10000[, [1000, 3000[, < 1000 | GDPH1, GDPH2, GDPH3, GDPH4 |
| **Unemployment rate** | [0, 10[, [10, 20[, ≥ 20 | UNEM1, UNEM2, UNEM3 |
| **Inflation rate** | [0, 10[, [10, 50[, [50, 100[, ≥ 100 | INFL1, INFL2, INFL3, INFL4 |
| **IHD** | 1, 2, 3 | IHD1, IHD2, IHD3 |

Table 3: *The qualitative variables for the POP_96 database.*

There are 8 qualitative variables and 31 modalities.

In Fig. 4, we represent the result of KMCA applied to these data. We can observe that the “good” modalities (level ended by 1), characteristic of developed countries, are displayed together in the bottom left-hand corner, followed by intermediate modalities (level 2), in the top left-hand corner. The bad modalities are displayed in the top right-hand corner, and the very bad ones in the bottom right-hand corner. There are 3 empty classes which separate the best modalities from the worst ones.

| | | | | | |
|---|---|---|---|---|---|
| ANPR2<br>UNEM2 | | MORT2<br>ANALR2 | SCHO2 | | MORT3<br>ANALR3 |
| | GDPH2 | | GDPH3 | INFL2<br>IHD2 | ANPR3 |
| UNEM1 | | INFL3 | INFL4 | | UNEM3 |
| IHD1 | INFL1 | | ANALR4 | ANPR4 | |
| ANPR1<br>ANALR1<br>SCHO1 | | | MORT4 | IHD3 | MORT5<br>ANALR5 |
| | MORT1<br>GDPH1 | | GDPH4 | FSCHO3 | |

Fig. 4: *The repartition of the 31 modalities on the Kohonen map by KMCA, 500 iterations.*

If we examine the projection (Fig. 5) on the first two axes by a Multiple Component Analysis, we see the same rough display, with a first axis which opposes the best modalities (on the left) to the worst ones (on the right). The projection on axes 1 and 5 (Fig. 6) allows a better legibility, but with a loss of information (only 30% of explained inertia). To keep 80% of the total inertia, it is necessary to consider the first 10 axes.

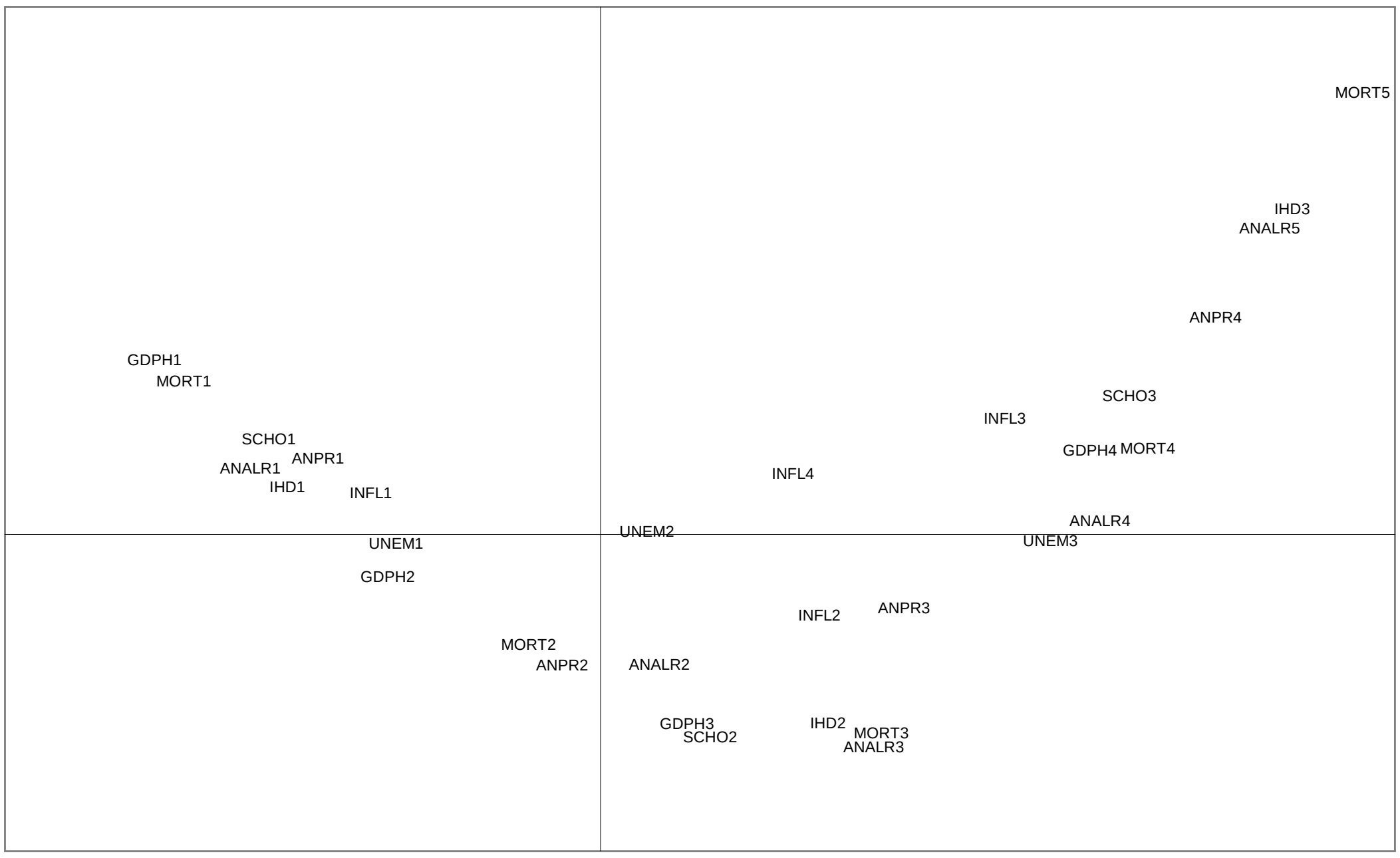

Fig 5: *The MCA representation, axes 1 (24%), and 2 (14%), 38% of explained inertia.*

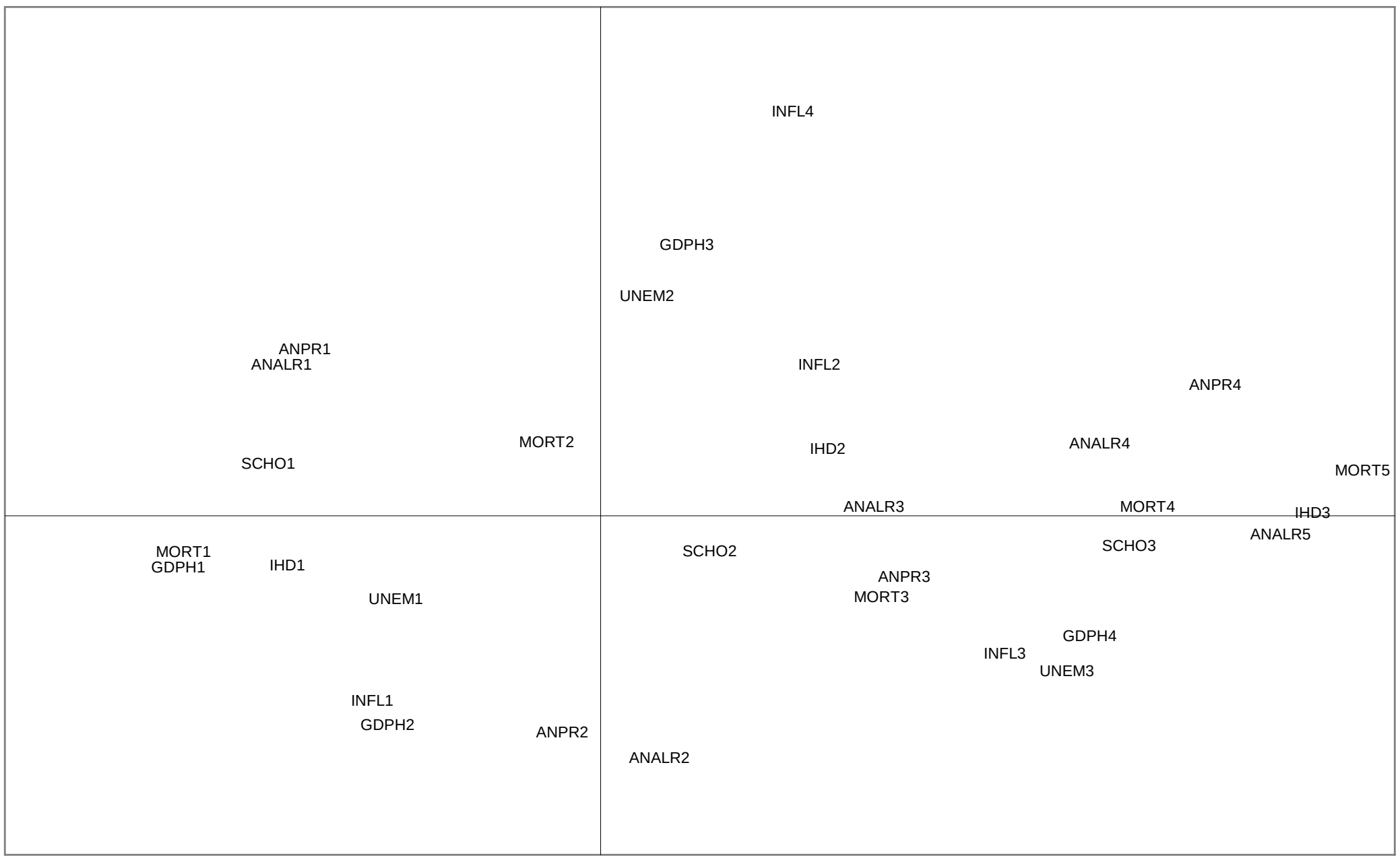

Fig 6: *The MCA representation, axes 1 (24%), and 5 (6%), 30% of explained inertia.*

As the modalities are classified according to some (irregular) scale, the progression from level 1 (the best) to the worst (level 3, or 4, or 5, depending on the modality) is clearly visible. But the clusters are not clear. The examination of further axes could give some contradictory conclusions. For example, modalities INFL3 and INFL4 appear to be close in Fig 5, and far in Fig. 6.

It is also possible to reduce the number of classes as in section 2.3, and to make this two-level classification visible, we allot to each « macro-class » some color or grey level. See Fig.7, a representation of the « micro-classes » grouped together to constitute 6 « macro-classes ».

| | | | | | |
|---|---|---|---|---|---|
| ANPR2<br>UNEM2 | | MORT2<br>ANALR2 | SCHO2 | | MORT3<br>ANALR3 |
| | GDPH2 | | GDPH3 | INFL2<br>IHD2 | ANPR3 |
| UNEM1 | | INFL3 | INFL4 | | UNEM3 |
| IHD1 | INFL1 | | ANALR4 | ANPR4 | |
| ANPR1<br>ANALR1<br>SCHO1 | | | MORT4 | IHD3 | MORT5<br>ANALR5 |
| | MORT1<br>GDPH1 | | GDPH4 | FSCHO3 | |

Fig. 7: *Macro-classes grouping the modalities into 6 easy-to-describe classes.*

## 7. Analyzing qualitative variables while keeping the individuals : KMCA_ind

Until now we have been interested only in the relations between modalities. But it can be more interesting and more valuable to cluster the individuals and the modalities that describe them at the same time. In this case, it is obvious that we need to deal with the Complete Disjunctive Table, in order to know the individual answers.

An easy way to simultaneously represent the individuals together with the modalities is to build a Kohonen map with the individuals using the algorithm SOM applied to the Corrected Complete Disjunctive Table, and then to project the modalities as supplementary data, with the suitable scaling.

Therefore, table $D$ is corrected into $D^c$ (following section 4.2), the SOM algorithm is trained with the rows of this corrected table. Then each modality $j$ is represented by an $M$-vector, which is the **mean vector** of all the individuals who share this modality. Its coordinates are:

$$\frac{b_{jl}}{b_j \sqrt{b_l}\sqrt{K}}, \quad \text{for} \qquad l = 1, \cdots, M,$$

with the notations defined in section 3. Each mean vector is assigned into the Kohonen class of its nearest code vector. This method is denoted by KMCA_ind and provides a simultaneous representation of individuals and modalities.

More iterations are necessary for the learning than for a simple classification of the modalities with KMCA. This is easily understandable since we have to first classify the individuals (which are generally more numerous than the modalities), but also to accurately compute the code vectors which are used as prototypes to assign the modalities considered as supplementary data.

The visualization is not always useful, especially when there are too many individuals. But the map provides visual associations between modalities or groups of modalities and subsets of individuals. See Ibbou, Cottrell (1995), Cottrell, de Bodt and Henrion (1996) for details and applications.

See in Fig. 8 the simultaneous representation of the 96 countries and the 31 modalities. Applying an Ascending Hierarchical, we have grouped them into 7 macro-classes.

| | | | | | |
|---|---|---|---|---|---|
| INFL4<br>Moldavia<br>Romania<br>Russia<br>Ukraine | Bulgaria<br>Poland | GDPH3<br>Costa Rica<br>Ecuador<br>Jamaica<br>Lebanon | Colombia<br>Fiji<br>Panama<br>Peru<br>Thailand | | MORT5<br>Afghanistan<br>Angola<br>Haiti<br>Mozambique<br>Pakistan<br>Yemen |
| Brazil | ANPR2<br>MORT2 | Croatia<br>Venezuela | ANALR2<br>UNEM2 | Ghana<br>Mauritania<br>Sudan | ANPR4<br>ANALR5<br>IHD3 |
| GPDH2<br>Chili<br>Cyprus<br>S. Korea | Argentina<br>Bahrain<br>Malaysia<br>Malta<br>Mexico | | INFL3<br>Macedonia<br>Mongolia | SCHO3<br>GPDH4<br>UNEM3<br>Laos | MORT4<br>ANALR4<br>Cameroon<br>Comoros<br>Ivory Coast<br>Nigeria |
| Greece<br>Hungary<br>Slovenia<br>Uruguay | UNEM1<br>IHD1<br>Portugal | China<br>Philippines<br>Yugoslavia | Albania<br>Indonesia<br>Sri Lanka | INFL2<br>Guyana<br>Vietnam | Bolivia<br>Kenya<br>Nicaragua |
| ANPR1<br>ANALR1<br>SCHO1<br>R. Czech | Germany, Sweden<br>Australia, Israel<br>USA, Iceland<br>Japan, Norway<br>New-Zealand<br>Netherlands<br>United Kingdom<br>Switzerland | INFL1<br>U Arab Emirates | SCHO2<br>Egypt | ANPR3<br>IHD2<br>Morocco<br>Paraguay | ANALR3<br>El Salvador<br>Swaziland |
| Belgium<br>Canada<br>Denmark<br>Spain, Italy<br>Finland<br>France<br>Ireland | MORT1<br>GDPH1<br>Luxemburg | Singapore | Saudi Arabia<br>Syria | MORT3<br>Algeria | South Africa<br>Iran<br>Namibia<br>Tunisia<br>Turkey<br>Zimbabwe |

Fig. 8: *The simultaneous representation of the 96 countries and the 31 modalities, 20 000 iterations.*

We observe that the proximity between countries and modalities is perfectly coherent, rich countries are in the bottom left-hand corner, with almost all the “good” modalities (level 1), very poor ones are in the top right-hand corner (level 3, 4, 5), etc.

This method achieves a nice simultaneous representation of modalities and individuals, but breaks the symmetry between individuals and modalities. Conversely, the algorithm presented in the next section keeps this symmetry and is directly inspired by classical methods. It was introduced in Cottrell et Letrémy (2003).

## 8. A new algorithm for a simultaneous analysis of individuals together the modalities: KDISJ

To keep together modalities and individuals in a more balanced way, in the same manner as for classical MCA, we define a new algorithm: KDISJ. In fact it is an extension of the algorithm KORRESP, (Cottrell et al., 1993) that has been introduced to analyze contingency tables which cross two qualitative variables. Algorithm KDISJ was first defined in (Cottrell, Letrémy, 2003). and an extended version is accepted for publication (Cottrell, Letrémy, 2004).

The Complete Disjunctive Table is corrected as shown in section 4.2. We then choose a Kohonen network, and associate with each unit $u$ a code vector $C_u$ that is comprised of $(M + N)$ components, with the first $M$ components evolving in the space for individuals (represented by the rows of $D^c$) and the $N$ final components in the space for modalities (represented by the columns of $D^c$ ).

We put down

$$C_u = (C_M, C_N)_u = (C_{M,u}, C_{N,u})$$

to highlight the structure of the code-vector $C_u$. The Kohonen algorithm lends itself to a double learning process. At each step, we alternatively draw a $D^c$ row (i.e. an individual *i*), or a $D^c$ column (i.e. a modality *j*).

When we draw an individual *i*, we associate the modality *j*(*i*) defined by

$$j(i) = Arg \max_j d^c_{ij} = Arg \max_j \frac{d_{ij}}{\sqrt{K d_{\cdot j}}}$$

that maximizes the coefficient $d^c_{ij}$, i.e. the rarest modality out of all of the corresponding ones in the total population. This modality is the most characteristic for this individual. In case of *ex-aequo*, we do as usual in this situation, by randomly drawing a modality among the candidates. We then create an extended individual vector $X = (i, j(i)) = (X_M, X_N)$, of dimension $(M + N)$. See Fig. 9. Subsequently, we look for the closest of all the code vectors, in terms of the Euclidean distance restricted to the first $M$ components. Let us denote by $u_0$ the winning unit. Next, we move the code vectors of the unit $u_0$ and its neighbors closer to the extended vector $X = (i, j(i))$, as per the customary Kohonen law. Let us write down the formal definition :

$$\begin{cases} u_0 = Arg \min_u \left\| X_M - C_{M,u} \right\| \\ C_u^{new} = C_u^{old} + \varepsilon \ \ \sigma(u, u_0) \ \ (X - C_u^{old}) \end{cases},$$

where $\varepsilon$ is the adaptation parameter (positive, decreasing with time), and $\sigma$ is the neighborhood function, defined as usual by $\sigma(u, u_0) = 1$ if $u$ and $u_0$ are neighbour in the Kohonen network, and = 0 if not. The adaptation parameter and the radius of the neighborhood vary as defined in the note 1, section 5.

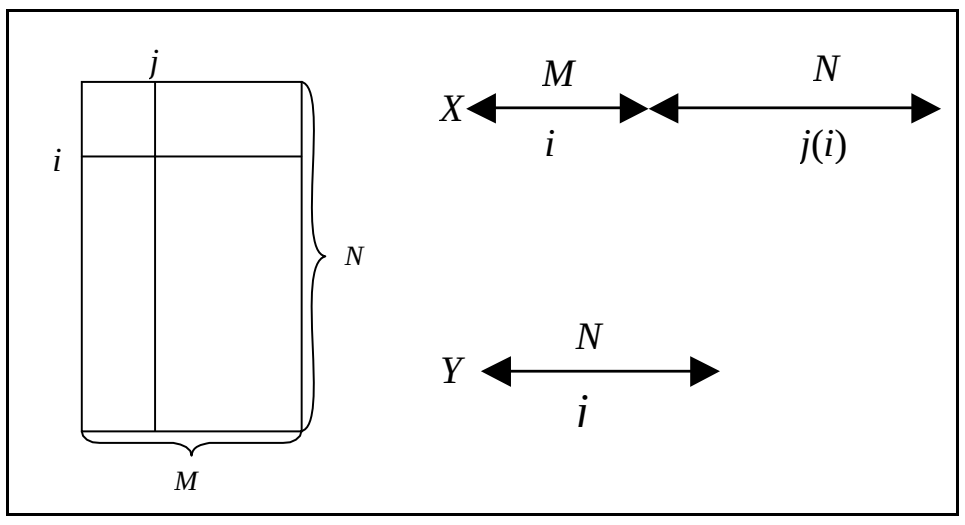

Fig 9 *: The matrix $D^c$, vectors X and Y.*

When we draw a modality *j* with dimension *N* (a column of $D^c$), we do not associate an individual with it. Indeed, by construction, there are many equally-placed individuals, and this would be an arbitrary choice. We then seek the code vector that is the closest, in terms of the Euclidean distance restricted to the last *N* components. Let $v_0$ be the winning unit. We then move the last *N* components of the winning code vectors associated to $v_0$ and its neighbors closer to the corresponding components of the modality vector *j*, without modifying the *M* first components. For the sake of simplicity let us denote by *Y* (see Fig. 9) the *N*-column vector corresponding to modality *j*. This step can be written :

$$\begin{cases} v_0 = Arg \min_u \left\| Y - C_{N,u} \right\| \\ C_{N,u}^{new} = C_{N,u}^{old} + \varepsilon \quad \sigma(u, v_0) \quad (Y - C_{N,u}^{old}) \end{cases}$$

while the first *M* components are not modified. See in Fig 10, the result of KDISJ applied to the POP_96 database.

| | | | | | |
|---|---|---|---|---|---|
| SCHO2<br>IHD2<br>Algeria<br>Syria | Saudi Arabia<br>Egypt<br>Indonesia | Brazil<br>Mexico | Argentina<br>Chili<br>Cyprus<br>S. Korea | INFL1<br>Australia<br>Canada<br>USA | GDPH1<br>Belgium<br>Denmark<br>Finland<br>France<br>Ireland<br>Italy |
| ANALR3<br>South Africa<br>Iran<br>Namibia<br>El Salvador<br>Zimbabwe | MORT3<br>Guyana<br>Morocco<br>Paraguay<br>Tunisia<br>Turkey | | GDPH2 | IHD1<br>Israel | MORT1<br>Germany, U Kingdom<br>Iceland, Japan<br>Lux, Singapore<br>Norway, Sweden<br>New Zealand<br>Netherlands, Spain<br>Switzerland |
| Kenya<br>Nicaragua | Swaziland | UNEM2<br>U Arab Emirates<br>Malaysia | Malta<br>Portugal | UNEM2<br>Greece<br>Hungary<br>Slovenia<br>Uruguay | ANPR1<br>ANALR1<br>SCHO1<br>R. Czech |
| ANALR4<br>Comoros<br>Ivory Coast | MORT4<br>Cameroon<br>Nigeria | ANPR3<br>Bolivia | INFL2<br>Bahrain<br>Philippines | Poland | INFL4<br>Croatia<br>Moldavia<br>Romania<br>Russia<br>Ukraine |
| ANPR4<br>IHD3<br>Ghana | SCHO3<br>GDPH4<br>Laos<br>Mauritania<br>Sudan | UNEM3<br>Vietnam<br>Yugoslavia | | MORT2<br>ANALR2<br>Bulgaria<br>Ecuador<br>Jamaica<br>Lebanon, Peru | GDPH3<br>Costa Rica |
| MORT5<br>ANALR5<br>Afghanistan<br>Angola<br>Pakistan<br>Yemen | Haiti<br>Mozambique | INFL3<br>Macedonia<br>Mongolia | Albania<br>China<br>Sri Lanka | Colombia<br>Panama | ANPR2<br>Fiji<br>Thailand<br>Venezuela |

Fig. 10: *Kohonen map with simultaneous representation of modalities and countries. The 36 micro-classes are clustered into 7 macro-classes, 2 400 iterations.*

We can observe that the simultaneous positions of countries and modalities are meaningful, and that the macro-classes are easy to interpret. It is possible to control the good position of the modalities with respect to the individuals, by computing the deviations inside each Kohonen class. The deviation for a modality $l$ (shared by $b_l$ individuals) and for a class $k$ (with $n_k$ individuals) can be calculated as the difference between the number of individuals who possess this modality and belong to the class $k$ and the “theoretical ” number $b_l\, n_k / N$ which would correspond to a distribution of the modality $l$ in the class $k$ that matches its distribution throughout the total population. These deviations are all positive, and it means that there is attraction between the modalities and the individuals who belong to the same class.

In so doing, we are carrying out a classical Kohonen clustering of individuals, plus a clustering of modalities, all the while maintaining the associations of both objects. After convergence, the individuals and the modalities belong to the same Kohonen classes. “Neighboring” individuals or modalities are inside the same class or neighboring classes. We call the algorithm we have just defined KDISJ. Its computing time is short, the number of iterations is about 20 times the total number of individuals and modalities, much less than for the KCMA_ind algorithm.

We can compare Fig. 10 with the representations that we get using a classical MCA, Fig. 11 and 12. The rough characteristics are the same, but in Fig. 11 and 12, the clusters are not clear, there are some contradictions from one projection to another (see for example the locations of the modalities INFL3 and INFL4). The Kohonen map is much more useful to visualize.

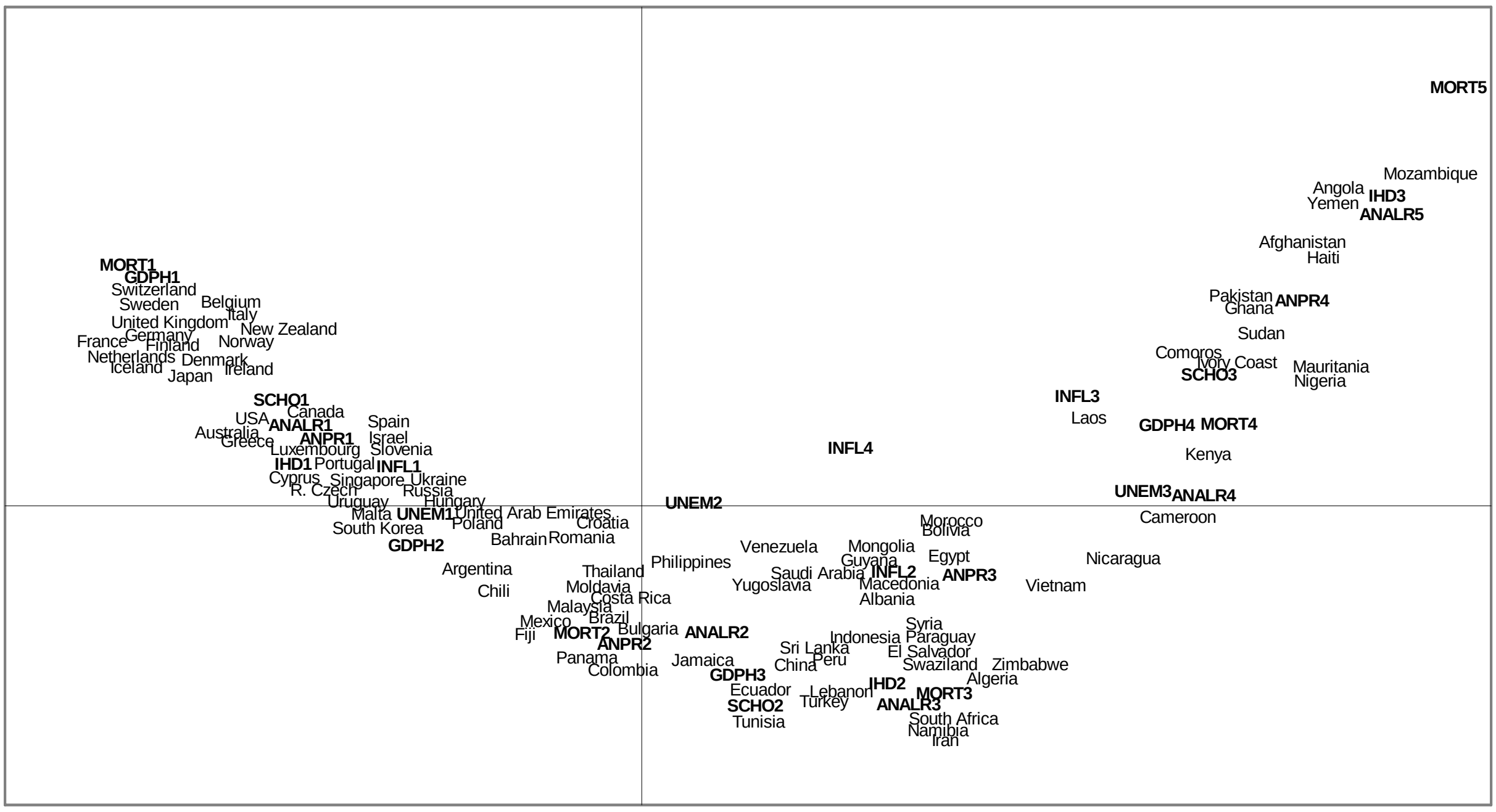

Fig. 11: *The MCA representation (modalities and individuals),*

*axes 1 (24%), and 2 (14%), 38% of explained inertia.*

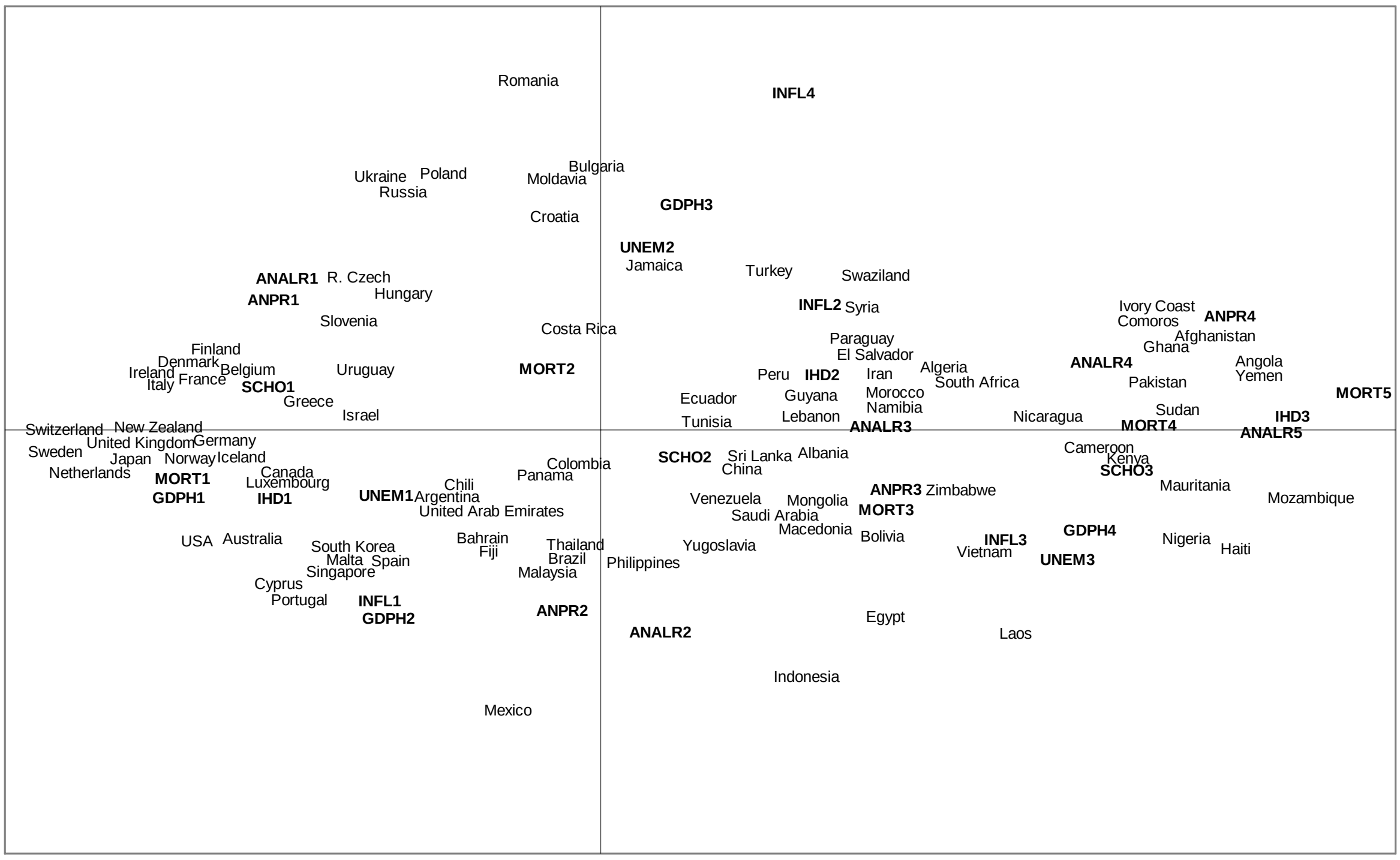

Fig. 12: *The MCA representation (modalities and individuals),*

*axes 1 (24%), and 5 (6%), 30% of explained inertia.*

## 9. Example II: a toy example, marriages

We consider 270 couples, the husband and the wife being classified into 6 professional groups : farmer, craftsman, manager, intermediate occupation, clerk, worker (numbered from 1 to 6).

These data are particularly simple since the contingency table has a very dominant principal diagonal. One knows that most marriages are concluded within the same professional category. See the data summarized in table 4. The complete disjunctive table is not displayed, but it is very easy to compute since there are only two variables (professional group of the husband, and of the wife). Among the 36 possible combinations for the couples, only 12 are present. So this example, while being a real-world data example, can be considered as a toy example.

| | FFARM | FCRAF | FMANA | FINTO | FCLER | FWORK | Total |
|---|---|---|---|---|---|---|---|
| **MFARM** | 16 | 0 | 0 | 0 | 0 | 0 | |
| **MCRAF** | 0 | 15 | 0 | 0 | 12 | 0 | 37 |
| **MMANA** | 0 | 0 | 13 | 15 | 12 | 0 | 40 |
| **MINTO** | 0 | 0 | 0 | 25 | 35 | 0 | 60 |
| **MCLER** | 0 | 0 | 0 | 0 | 25 | 0 | 25 |
| **MWORK** | 0 | 0 | 0 | 10 | 60 | 32 | 102 |
| Total | 16 | 15 | 13 | 50 | 144 | 32 | 270 |

Table 4: *Contingency Table for the married couples, (extracted for French INSEE statistics, in 1990). The rows are for the husbands, the columns for the wives.*

In what follows the couples are denoted by $(i, j)$ where $i = 1, ..., 6$, $j = 1, ..., 6$ correspond to the different groups as indicated above.

Let's first consider the results of a Factorial Correspondence Analysis. 5 axes are needed to keep 80% of the total inertia. The best projection, on axes 1 and 2 in Fig. 13. is too schematic. Axis 1 opposes the farmers to all other couples but distorts the representation.

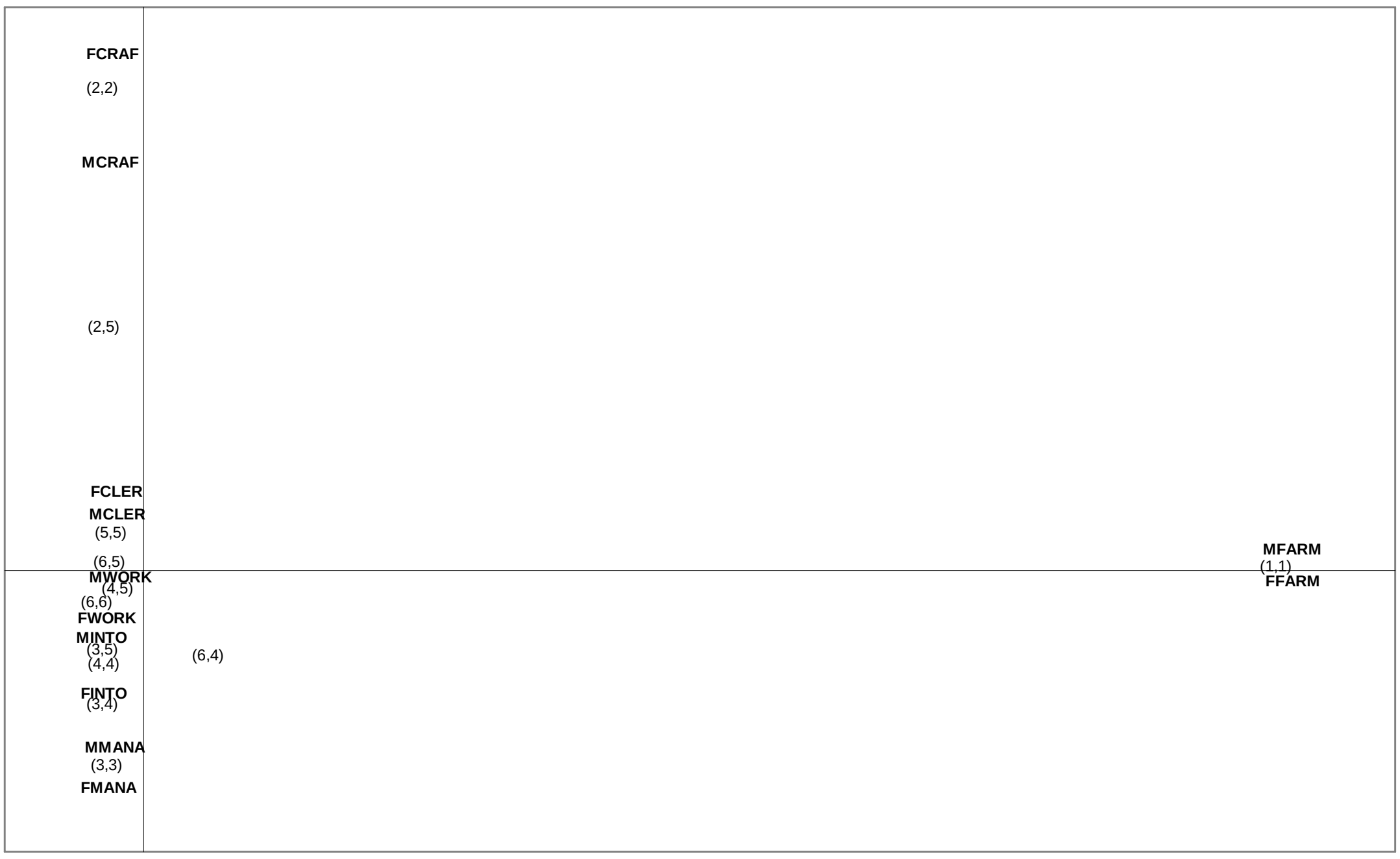

Fig. 13: *The MCA representation (modalities and individuals),*

*axes 1 (20%), and 2 (17%), 37% of explained inertia.*

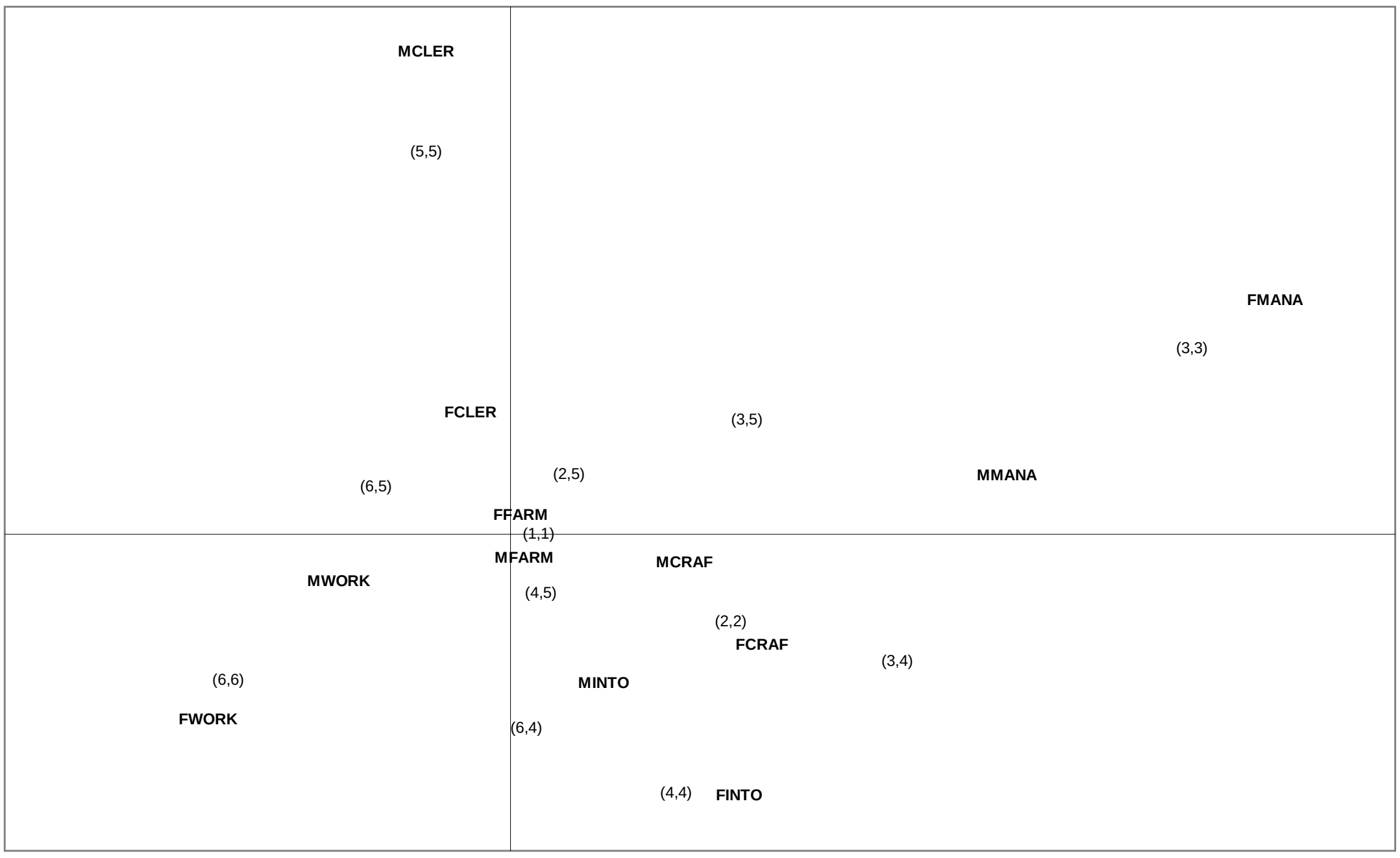

Fig. 14: *The MCA representation (modalities and individuals),*

*axes 3 (16%), and 5 (13%), 29% of explained inertia.*

On axes 3 and 5, the representation is better, even if the percentage of explained inertia is smaller. In both cases, each kind of couples is precisely put halfway between the modalities of each member of the couple.

Only 12 points are visible for individuals, since couples corresponding to the same professional groups for the husband and the wife are identical (there are no other variables).

Let us apply the KCMA method to these very simple data. We get the SOM map in Fig. 15.

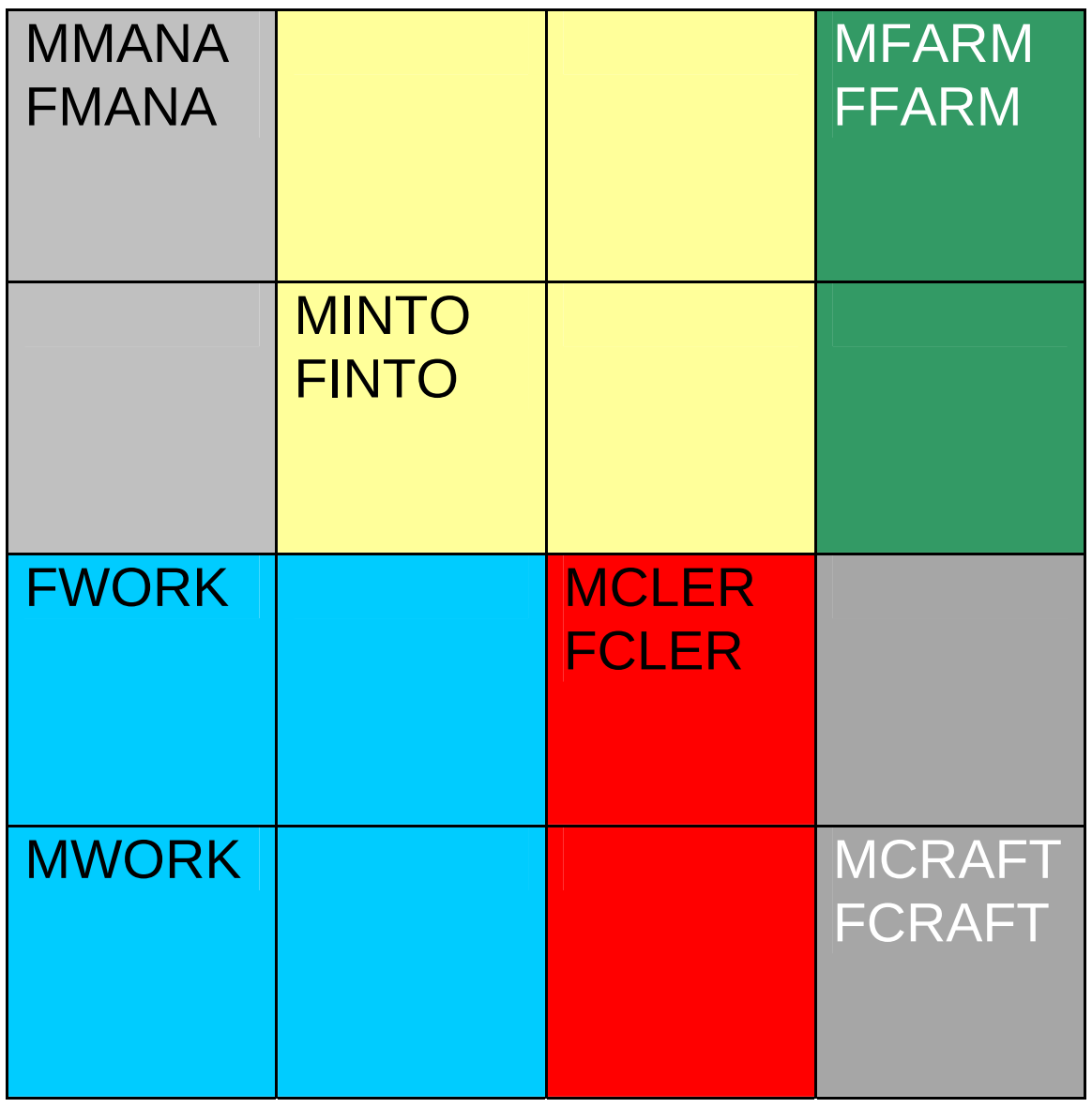

Fig. 15: *Kohonen map with representation of modalities. The 16 micro-classes are clustered into 6 macro-classes which gather only identical modalities, 200 iterations.*

In what follows (Fig. 16), we used KCMA_ind algorithm. We indicate the number of each kind of couples present in each class.

| | | | |
|---|---|---|---|
| MFARM<br>FFARM<br><br>16 (1,1) | | FMANA<br><br>13 (3,3) | MMANA |
| | MCLER<br><br>25 (5,5) | | 15 (3,4)<br>12 (3,5) |
| MCRAFT<br>FCRAFT<br><br>15 (2,2)<br>12 (2,5) | FCLER | 60 (6,5) | MWORK |
| | MINTO<br><br>25 (4,4)<br>35 (4,5) | FINTO<br><br>10 (6,4) | FWORK<br><br>32 (6,6) |

Fig. 16: KCMA_ind : *Kohonen map with simultaneous representation of modalities and individuals. The number of couples of each type is indicated. 16 (1, 1) means that there are 16 couples where the husband is a farmer and the woman too, 10 000 iterations.*

In Fig. 17, we represent the results obtained with a KDISJ algorithm. The results are very similar. In both Kohonen maps, each type of couples is situated in the same class as the professional groups of both members of the couple, or between the corresponding modalities.

| | | | |
|---|---|---|---|
| MFARM<br>FFARM<br><br>16 (1,1) | | FINTO<br><br>25 (4,4)<br>10 (6,4) | MMANA<br><br>15 (3,4) |
| | MINTO<br><br>35 (4,5) | | FMANA<br><br>13 (3,3)<br>12 (3,5) |
| MCRAFT<br>FCRAFT<br><br>15 (2,2)<br>12 (2,5) | FCLER | MWORK<br><br>60 (6,5) | |
| | MCLER<br><br>25 (5,5) | | FWORK<br><br>32 (6,6) |

Fig. 17: KDISJ : *Kohonen map with simultaneous representation of modalities and individuals. The number of couples of each sort is indicated: 16 (1, 1) means that there are 16 couples where the man is a farmer and the woman too, 5000 iterations.*

In this toy example we can therefore conclude that the results are quite good, since they give the same information as the linear projections, with the advantage that only one map is sufficient to summarize the structure of the data.

## 10. Example III: temporary agency contracts

In this section, we present another example extracted from a large study of the INSEE's 1998-99 Timetable survey. The complete report (Letrémy, Macaire et. al., 2002) is an attempt at determining which working patterns have a specific effect and to question the homogeneity of the form of employment. Specifically, they study the "Non Standard Contracts" that is i) the various temporary contracts, or fixed-term contracts, ii) the various part-time contracts, fixed-term or unlimited-term contracts, iii) the temporary agency work contracts. The paper analyzes which specific time constraints are supported by non-standard employment contracts. Do they always imply a harder situation for employees, compared to standard employment contracts.

In the survey, the employees had to answer some questions about their "Working times". There are varied issues : working time durations, schedules and calendars, work rhythms, variability, flexibility and predictability of all these time dimensions, possible choices for employees, etc.

An initial study entitled “Working times in particular forms of employment: the specific case of part-time work” (Letrémy, Cottrell, 2003) covered 14 of the questions that the questionnaire had asked, representing 39 response modalities and 827 part-time workers.

In this section, we study the employees who have temporary agency work contracts, there are 115 employees having this kind of contract. We present an application of KMCA and KDISJ algorithms that are used to classify the modalities of the survey as well as the individuals.

Table 5 lists the variables and response modalities that were included in this study. There are 25 modalities.

| **Heading** | **Name** | **Response modalities** |
|---|---|---|
| Sex | **Sex 1 2** | Man, Woman |
| Age | **Age 1, 2, 3, 4** | <25, [25, 40[, [40,50[, ≥50 |
| Daily work schedules | **Dsch 1, 2, 3** | Identical, as-Posted, Variable |
| Number of days worked in a week | **Dwk 1, 2** | Identical, Variable |
| Night work | **Night 1, 2** | No, Yes |
| Saturday work | **Sat 1, 2** | No, Yes |
| Sunday work | **Sun 1, 2** | No, Yes |
| Ability to go on leave | **Leav 1, 2, 3** | Yes no problem, yes under conditions, no |
| Awareness of next week schedule | **Nextw 1, 2** | Yes, no |
| Possibility of carrying over credit hours | **Car 0, 1, 2** | No point, yes, no |

Table 5: *Variables that were used in the individual survey.*

In Fig. 18, the modalities are displayed on a Kohonen map, after being classified by a KMCA algorithm.

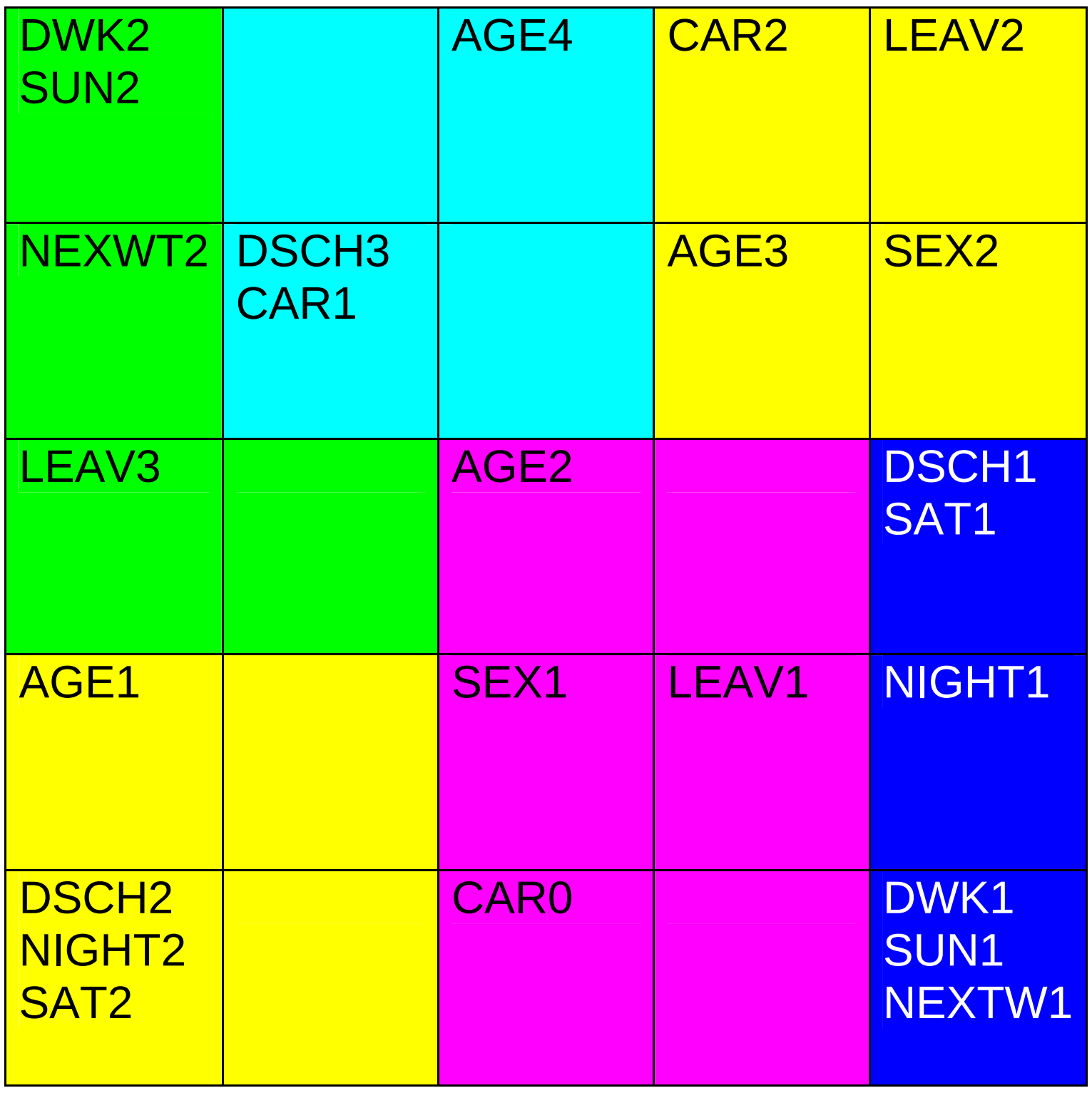

Fig 18: *The modalities are displayed; they are grouped into 6 clusters, 500 iterations.*

The clusters are clearly identifiable : the best work conditions are in the bottom right-hand corner (level number 1), the youngest people (AGE1) are in the bottom left-hand corner associated with bad conditions (they work Saturdays and nights, etc.)

One can see the same associations on the MCA representation, see Fig. 19 for axes 1 and 2.

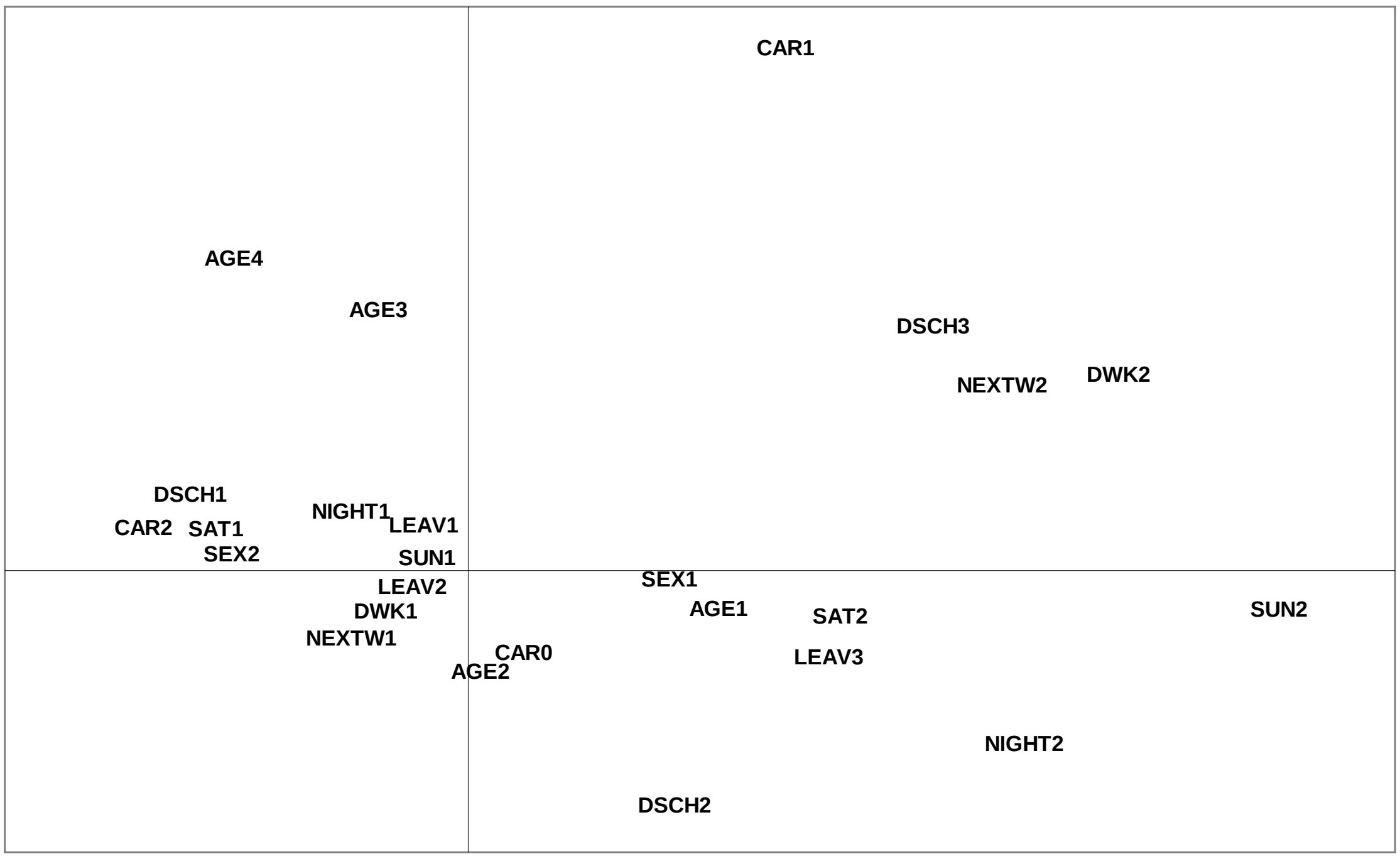

Fig 19: *The MCA representation, axes 1 (19%), and 2 (11%), 30% of explained inertia.*

Now, we can simultaneously classify the modalities together with the individuals by using the KDISJ algorithm. See Fig. 20.

| | | | | |
|---|---|---|---|---|
| DWK2<br>6 ind | 3 ind | DSCH3<br>7 ind | CAR1<br>1 ind | AGE3<br>9 ind |
| 2 ind | NEXTW2<br>2 ind | LEAV3<br>3 ind | 2 ind | 4 ind |
| SUN2<br>6 ind | 2 ind | 4 ind | 5 ind | LEAV2<br>6 ind |
| 4 ind | SEX1<br>DSCH2<br>CAR0<br>4 ind | AGE2<br>DWK1<br>SUN1<br>NEXTW1<br>7 ind | SEX2<br>DSCH1<br>NIGHT1<br>SAT1<br>3 ind | CAR2<br>7 ind |
| NIGHT2<br>SAT2<br>6 ind | AGE1<br>9 ind | LEAV1<br>5 ind | 4 ind | AGE4<br>4 ind |

Fig. 20: KDISJ: *Simultaneous classification of 25 modalities and 115 individuals. Only the number of individuals is indicated in each class. 3000 iterations.*

The groups are easy to interpret, the good conditions of work are gathered, the bad ones as well.

As usual, the projection of both individuals and modalities on the first two axes of a MCA is not very clear and the visualization is poor. See Fig.21, only 30% of the total inertia is taken into account and 9 axes would be needed to obtain 80% of the total inertia.

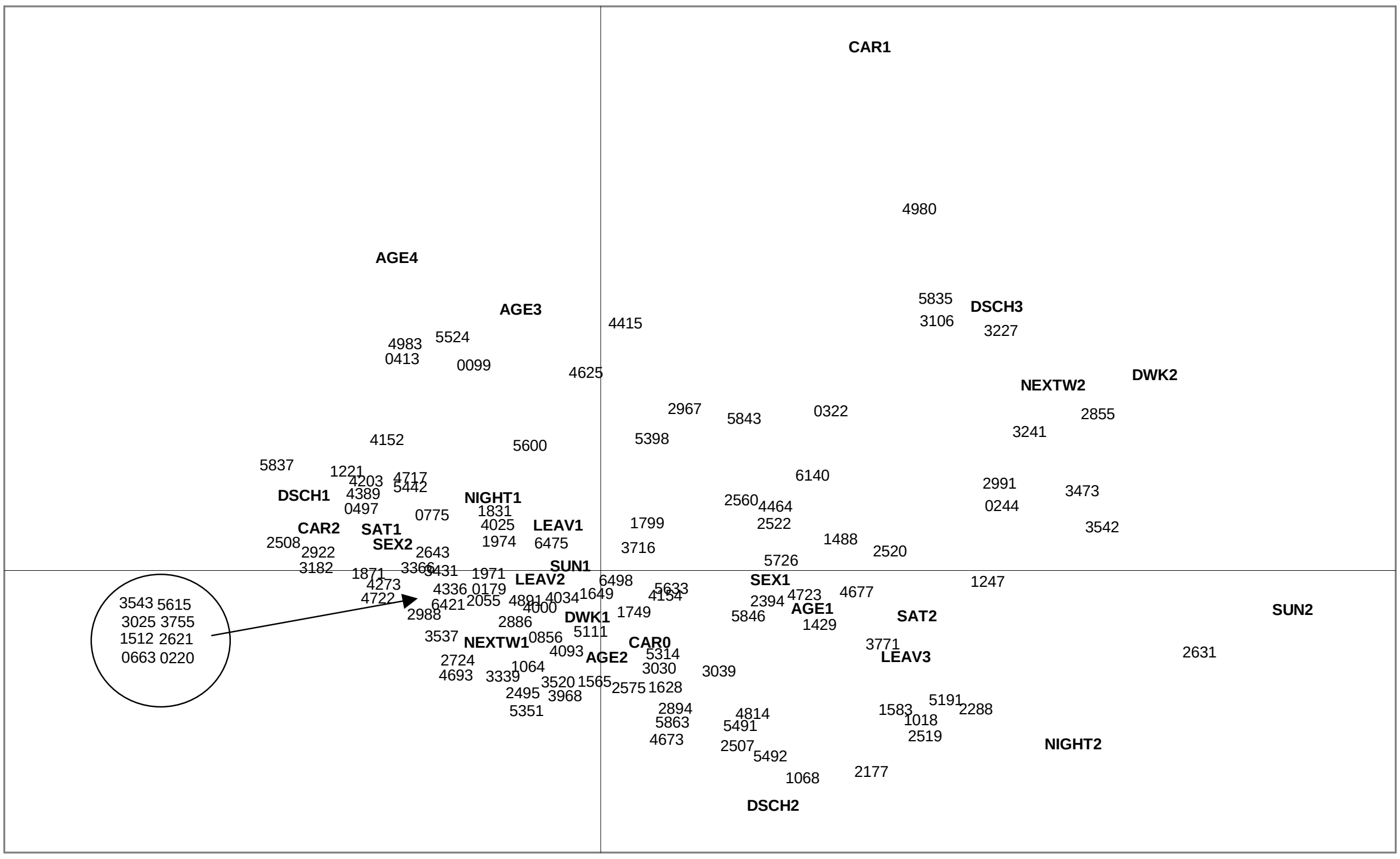

Fig. 21: *Simultaneous MCA representation of the 25 modalities and the 115 individuals, axes 1 (19%), and 2 (11%), 30% of explained inertia.*

## 11. Conclusion

We propose several methods to analyze multidimensional data, in particular when observations are described by qualitative variables, as a complement of classical linear and factorial methods. These methods are adaptations of the original Kohonen algorithm.

Let us summarize the relations between the classical factorial algorithms and the SOM-based ones in table 6.

<table>
<tr><th>SOM-based algorithm</th><th>Factorial method</th></tr>
<tr><td>SOM algorithm on the rows of matrix X</td><td>PCA, diagonalization on $X'X$.</td></tr>
<tr><td>KCMA (clustering the modalities): SOM algorithm on the rows of $B^c$</td><td>MCA, diagonalization on $B^{c}{}'B^c$.</td></tr>
<tr><td>KCMA_ind (clustering the individuals): SOM algorithm on the rows of $D^c$ and setting of the modalities,</td><td rowspan="2">MCA with individuals, diagonalization of $D^{c}{}'D^c$ and $D^cD^{c}{}'$.</td></tr>
<tr><td>KDISJ (coupled training with the rows -individuals- and the columns -modalities): SOM on the rows and the columns of $D^c$</td></tr>
</table>

Table 6: *Comparison between the SOM methods and the Factorial Methods.*

But in fact, for applications, it is necessary to combine different techniques. For example, in case of quantitative variables, it is often interesting to first reduce the dimension by applying a Principal Component Analysis, and by keeping a reduced number of coordinates.

On the other hand, if the observations are described with quantitative variables as well as qualitative ones, it is useful to build a classification of the observations restricted only the quantitative variables, using a Kohonen classification followed by an Ascending Hierarchical Algorithm, to define a new qualitative variable. It is added to the other qualitative variables and it is possible to apply a Multiple Correspondence Analysis or a KCMA to all the qualitative variables (the original variables and the class variable just defined). This technique leads to an easy description of the classes, and highlights the proximity between modalities.

If we are interested in the individuals only, the qualitative variables can be transformed into real-valued variables by a Multiple Correspondence analysis. In that case, all the axes are kept, and each observation is then described by its factorial coordinates. The database thus becomes numerical, and can be analyzed by any classical classification algorithm, or by a Kohonen algorithm.

In this paper, we do not give any example of one-dimensional Kohonen map. But when it is useful to establish a score of the data, the construction of a Kohonen string (of dimension 1) from the data or from the code vectors built from the data, straightforwardly gives a score by "ordering" the data.

It would be necessary to bear all these techniques in mind, together with classical techniques, to improve the performances of all of them, and consider them as very useful tools in data mining.